\documentclass[11pt,leqno,oneside]{article}

\usepackage{amsmath}
\usepackage{amsthm}
\usepackage{amstext}
\usepackage{amsopn}
\usepackage{amsbsy}

\newtheorem*{main-theorem}{Main Theorem}
\newtheorem{theorem}{Theorem}
\newtheorem{corollary}{Corollary}
\newtheorem{lemma}{Lemma}

\begin{document}

\title{Lower Bounds of the First Closed and Neumann Eigenvalues of Compact Manifolds with Positive Ricci Curvature
\thanks{2000 Mathematics Subject Classification Primary 58J50, 35P15; Secondary 53C21}}
\author{Jun LING}
\date{}

\maketitle

\begin{abstract}
We give new estimates on the lower bounds for the first closed or
Neumann eigenvalue for a compact manifold with positive Ricci
curvature in terms of the diameter and the lower bound of Ricci
curvature. The results improve the previous estimates.
\end{abstract}

\section{Introduction}\label{sec-intro}
For an n-dimensional closed Riemannian manifold whose Ricci
curvature has a positive lower bound $(n-1)K$ for some constant
$K>0$, A.~Lichnerowicz \cite{lich} gave a lower bound of the first
eigenvalue $\lambda$ of the Laplacian
\begin{equation}\label{nK-bound}
\lambda\geq nK.
\end{equation}
Under the same curvature assumption, Escobar \cite{es} proved that
if the compact manifold has a weakly convex boundary, the the
first non-zero Neumann eigenvalue of $M$ has the above lower bound
(\ref{nK-bound}) as well.

The above Lichnerowicz-type lower bound (\ref{nK-bound}) gives no
information when the above constant $K$ vanishes. In such case,
Li-Yau \cite{liy1} and Zhong-Yang \cite{zy} provided another lower
bound for the first non-zero eigenvalue of a closed manifold
\begin{equation}\label{lyzy}
\lambda\geq \frac{\pi^2}{d^2}.
\end{equation}

It is an interesting problem to find a unified lower bound of the
first closed or Neumann eigenvalue $\lambda$ in terms of the lower
bound $(n-1)K$ of the Ricci curvature and the diameter $d$ so that
the lower bound of the the first non-zero eigenvalue does not
vanish as $K$ vanishes. P.~Li conjectured a unified bound of the
first non-zero eigenvalue should be $\pi^2/d^2 + (n-1)K$. There
has been some work along this line, say \cite{yang} by D.~Yang,
and \cite{ling3} by the author that improved Yang's estimate for
the first Dirichlet eigenvalue in \cite{yang}. D.~Yang \cite{yang}
also give an estimate on the lower bound of closed or Neumann
eigenvalue
\begin{equation}\label{yang-bound}
\lambda \geq \frac{\pi^2}{d^2}+\frac{1}{4}(n-1)K.
\end{equation}
In this paper, we give some new estimates on the lower bound and
improve the above bound. Our main result is the following theorem.

\begin{theorem} \label{main-thm}
If $M$ is an n-dimensional, compact Riemannian manifold that has
an empty or none-empty boundary whose second fundamental form is
nonnegative with respect to the outward normal (i.e., weakly
convex). Suppose that Ricci curvature \textup{Ric}$(M)$ has a
lower bound $(n-1)K$ for some constant $K>0$, that is
\begin{equation}\label{ricci-bound}
\textup{Ric}(M)\geq (n-1)K>0.
\end{equation}
Then the first non-zero (closed or Neumann, which applies)
eigenvalue $\lambda$ of the Lalacian $\Delta$ on $M$ has the
following lower bound
\[
\lambda \geq \frac{\pi^2}{d^2}+\frac{3}{8}(n-1)K \qquad\textrm{for
}n=2
\]
and
\[
\lambda \geq \frac{\pi^2}{d^2}+\frac{31}{100}(n-1)K
\qquad\textrm{for }n\geq 3,
\]
where where $d$ is the diameter of $M$.
\end{theorem}

This estimate sharpens Yang's bound (\ref{yang-bound}). It is a
generalization of Li-Yau \cite{liy1} and Zhong-Yang \cite{zy}'s
result (\ref{lyzy}) for a closed manifold and it is better than
Lichnerowicz's bound (\ref{nK-bound}) if the manifold is
non-symmetric and has small diameter with respect to the positive
lower bound of the Ricci curvature.

In order to improve the known results, we need to construct
suitable test functions where detailed technical work is
essential. In Section \ref{sec-functions} we construct the test
function $\xi$. We explore the properties of the function $\xi$,
the Zhong-Yang function $\eta$, and the ratio function $\xi/\eta$.
Those properties are essential to the construction of the suitable
test functions. Because those functions are complicated
combinations of trigonometric and rational functions, the needed
properties such as monotonic and convex properties are hard to
prove. In the past, though we know that many nice properties might
be true, only a few of them could be proven strictly in
mathematics by the canonical calculus method and therefore be used
in strict mathematics proof. We are able to prove those properties
effectively now by studying the differential equations those
functions satisfied and using the Maximum Principle. Since the
constructions and proofs in that part are quite technical by
nature, we put them in the last section. Readers may refer to that
section when in need. We derive several preliminary estimates in
the next section and prove our result in Section \ref{sec-proofs}.

\section{Preliminary Estimates}\label{sec-pre-es}
The first preliminary estimate is due to Lichnerowicz and Escobar.
For the completeness and consistency, we use gradient estimate in
\cite{li}-\cite{liy1} and \cite{sy} to derive the two estimates.
\begin{lemma}               \label{nK-lemma}
Let $\lambda$ be the first non-zero (closed or Neumann, which
applies) eigenvalue under the conditions in Theorem
\ref{main-thm}. Then (\ref{nK-bound}) holds.
\end{lemma}
\begin{proof}
Let $u$ be a normalized eigenfunction of the first non-zero
(closed or Neumann, which applies) eigenvalue $\lambda$ such that
\[ \sup_{M}u=1,
\quad \inf_{M}u=-k, \quad \textrm{and}\quad 0<k\leq 1,
\]
and define a function $v$ by
\begin{equation}                                \label{v-n-def}
v=[ u - (1-k)/2 ]/ [ (1+k)/2 ].
\end{equation}

Then
\begin{equation}                                     \label{v-n-con}
\max v=1\quad \textrm{and} \quad \min v=-1.
\end{equation}
 The function $v$  satisfies the following equation
\begin{equation}                                       \label{v-n-eq}
\Delta v=-\lambda (v+a) \quad \textrm{in }M,
\end{equation}
where
\begin{equation}                                        \label{a-def}
a=(1-k)/(1+k).
\end{equation}
Note that $0\leq a<1$. If $M$ has non-empty boundary $\partial M$,
then $v$ satisfies Neumann condition on the boundary,
\begin{equation}                                         \label{v-n-boundary}
\frac{\partial v}{\partial \nu}=0  \qquad \textrm{on }\partial M,
\end{equation}
where $\nu$ is the the outward normal of $\partial M$.

Take an local orthonormal frame $\{e_1, \dots, e_n\}$ about
$x_0\in M$. At $x_0$
\[
\nabla_{e_j}(|\nabla v|^2)(x_0)= \sum_{i=1}^n2v_iv_{ij}
\]
and
\begin{eqnarray}
\Delta (|\nabla v|^2)(x_0)&=&2\sum_{i,j=1}^nv_{ij}v_{ij}+2\sum_{i,j=1}^nv_{i}v_{ijj}\nonumber\\
      &=&2\sum_{i,j=1}^nv_{ij}v_{ij}+2\sum_{i,j=1}^nv_{i}v_{jji}+2\sum_{i,j=1}^n\textup{R}_{ij}v_iv_j\nonumber\\
      &=&2\sum_{i,j=1}^nv_{ij}v_{ij}+2\nabla v \nabla (\Delta v) + 2\textup{Ric}(\nabla v, \nabla v)\nonumber\\
      &\geq&2\sum_{i=1}^nv_{ii}^2+ 2\nabla v \nabla (\Delta v) + 2(n-1)K|\nabla
      v|^2\nonumber\\
      &\geq&\frac{2}{n}(\Delta v)^2-2\lambda|\nabla v |^2 + 2(n-1)K|\nabla v|^2.\nonumber
\end{eqnarray}
Thus at all point $x\in M$,
\begin{equation}                        \label{basic2.002}
\frac12\Delta (|\nabla v|^2)\geq \frac1n \lambda^2(v+a)^2 +
[(n-1)K-\lambda]|\nabla v|^2.
\end{equation}
On the other hand, after multiplying (\ref{v-n-eq}) by $v+a$ and
integrating the both sides over $M$. When $M$ has non-empty
boundary and $v$ satisfies Neumann condition (\ref{v-n-boundary}),
we have
\[
\int_{M}\lambda(v+a)^2\,dx=-\int_M(v+a)\Delta v\, dx
\]
\[=-\int_{\partial M}(v+a)\frac{\partial}{\partial
\nu}v\,ds+\int_{M}|\nabla v|^2\, dx=\int_{M}|\nabla v|^2\, dx.
\]
That the integral on the boundary vanishes is due to
(\ref{v-n-boundary}). Integrating (\ref{basic2.002}) over $M$ and
using the above equality, we get
\begin{equation}\label{basic2.003}
\frac12\int_{\partial M}\frac{\partial}{\partial \nu}(|\nabla
v|^2)\,dx \geq\int_M(nK-\lambda)\frac{n-1}{n}\lambda(v+a)^2\,dx.
\end{equation}
We want to show that $\frac{\partial}{\partial \nu}(|\nabla
v|^2)\leq 0$ on $\partial M$. Take any $x_0\in \partial M$. If
$\nabla v(x_0)=0$, then it is done. Assume now that $\nabla
v(x_0)\not= 0$. Choose an orthonormal frame $\{ e_1, \dots, e_n\}$
about $x_0$ such that $e_n|_{x_0}$ is the unit outward normal
vector to $\partial M$ at $x_0$. Let $(h_{ij})$ be the second
fundamental form of $\partial M$ with respect to the outward
normal $\nu$ to $\partial M$. Now at $x_0$,
\begin{eqnarray}
v_{in}&=&e_ie_nv-(\nabla _{e_i}e_n)v\nonumber\\
      &=&-(\nabla _{e_i}e_n)v\nonumber\\
      &=&-\sum_{j=1}^{n-1} h_{ij}v_j\nonumber
\end{eqnarray}
and
\begin{eqnarray}
\frac{\partial}{\partial \nu}(|\nabla v|^2)&=&e_n|\nabla v|^2=2\sum_{j=1}^nv_jv_{jn}\nonumber\\
      &=&2\sum_{j=1}^{n-1}v_jv_{jn}=-2\sum_{i, j=1}^{n-1} h_{ij}v_iv_j\nonumber\\
      &\leq& 0 \qquad \textrm{by the weak convexity of }\partial
      M.\label{der-grad-v}
\end{eqnarray}
Putting this into (\ref{basic2.003}), we get the Lichnerowicz-type
bound (\ref{nK-bound}) for the first non-zero Neumann eigenvalue.
We get the bound (\ref{nK-bound}) for the first non-zero closed
eigenvalue by a similar argument as the above when $M$ has no
boundary, just noticing that there are no boundary terms in such
case.
\end{proof}

\begin{lemma}       \label{pre-es-n}
Let $v$ be the same as in (\ref{v-n-def}). Then $v$ satisfies the
following
\begin{equation}                        \label{basic3-n}
\frac{\left |\nabla v\right |^2}{b^2-v^2} \leq\lambda(1+a),
\end{equation}
where $a$ is defined in (\ref{a-def}) and $b>1$ is an arbitrary
constant.
\end{lemma}

\begin{proof}
 Consider the function
\begin{equation}                \label{p-of-x-def}
P(x)=|\nabla v|^2+Av^2,
\end{equation}
where $v$ is the function in (\ref{v-n-def}), and where $A=\lambda
(1+a)+\epsilon$ for small $\epsilon>0$. Function P must achieve
its maximum at some point $x_0\in M$.

We claim that
\begin{equation}\label{gra-p-eq-0}
\nabla P(x_0)=0.
\end{equation}

If $x_0\in M\backslash \partial M$, (\ref{gra-p-eq-0}) is
obviously true. Suppose that $x_0\in\partial M$. Choose a local
orthonormal frame $\{e_1, e_2,\cdots, e_{n}\}$ of $M$ about $x_0$
as in the proof of (\ref{der-grad-v}) such that $e_n$ is the unit
outward normal vector field near $x_0\in
\partial M$ and $\{e_1, e_2,\cdots, e_{n-1}\}|_{\partial M}$ is a
local frame of $\partial M$ about $x_0$. Note that
$\nabla_{e_n}e_i=0$ for $i\leq n-1$ and $v_n(x_0)=0$.

$P(x_0)$ is the maximum implies that
\begin{equation}\label{p-i-eq-0}
P_i(x_0)= 0\qquad \textrm{for }i\leq n-1
\end{equation}
and
\begin{equation}\label{p-n-geq-0}
P_n(x_0)\geq 0.
\end{equation}
Using (\ref{v-n-con})-(\ref{v-n-boundary}) in the following
arguments, then we have that at $x_0$,
\begin{eqnarray}
v_{in}&=&\sum_{i=1}^n e_ie_nv-\sum_{i=1}^n (\nabla _{e_i}e_n)v\nonumber\\
      &=&-\sum_{i=1}^n (\nabla _{e_i}e_n)v\nonumber\\
      &=&-\sum_{j=1}^{n-1} h_{ij}v_j\nonumber
\end{eqnarray}
and
\begin{eqnarray}\label{p-n-leq-0}
P_n&=&2\sum_{j=1}^nv_jv_{jn}+2Avv_n\nonumber\\
      &=&2\sum_{j=1}^{n-1}v_jv_{jn}=-2\sum_{i, j=1}^{n-1} h_{ij}v_iv_j\nonumber\\
      &\leq& 0 \qquad \textrm{by the convexity of }\partial M,
\end{eqnarray}
where $(h_{ij})$ is the second fundamental form of $\partial M$
with respect to the outward normal $e_n$.

Now (\ref{p-i-eq-0}), (\ref{p-n-geq-0}) and (\ref{p-n-leq-0})
imply that $P_n(x_0)=0$ and $\nabla P(x_0)=0$.

Therefore (\ref{gra-p-eq-0}) holds, no matter $x_0\not\in\partial
M$ or $x_0\in\partial M$. By (\ref{gra-p-eq-0}) and the Maximum
Principle, we have
\begin{equation}\label{max-prin}
\nabla P(x_0)=0 \qquad \textrm{and}\qquad \Delta P(x_0)\leq 0.
\end{equation}

There are two cases, either  $\nabla v(x_0)=0$ or $\nabla
v(x_0)\not=0$.

If $\nabla v(x_0)=0$, then
\[
|\nabla v(x)|^2+Av(x)^2=P(x)\leq P(x_0)\leq A.
\]
Let $\epsilon\rightarrow 0$ in the above inequality. Then
(\ref{basic3-n}) follows.

If $\nabla v (x_0)\not=0$, then we rotate the local orthonormal
frame about $x_0$ such that
\[
|v_1(x_0)|=|\nabla v(x_0)|\not=0\qquad \textrm{and}\qquad
v_i(x_0)=0,\quad i\geq2.
\]
From (\ref{max-prin}) we have at $x_0$,
\[
0=\frac12\nabla_i P=\sum_{i,j=1}^nv_jv_{ji}+\sum_{i=1}^nAvv_i,
\]
\begin{equation}                \label{d1}
v_{11}=-Av \qquad \textrm{and}\qquad v_{1i}=0\quad i\geq 2,
\end{equation}
and
\begin{eqnarray}
&0 &\geq \frac12\Delta P(x_0)= \sum_{i,j=1}^n \left(v_{ji}v_{ji}+v_{j}v_{jii}+Av_{i}v_{i}+Avv_{ii}\right)\nonumber\\
&{} &=\sum_{i,j=1}^n v_{ji}^2+\nabla v\nabla(\Delta v) +\textup{Ric}(\nabla v, \nabla v) +A|\nabla v|^2 +Av\Delta v\nonumber\\
&{} &\geq v_{11}^2+\nabla v\nabla(\Delta v) + (n-1)K|\nabla v|^2+A|\nabla v|^2 +Av\Delta v\nonumber\\
&{} &=(-Av)^2-\lambda |\nabla v|^2+ (n-1)K|\nabla v|^2+A|\nabla v|^2 -\lambda Av(v+a)\nonumber\\
&{} &=(A-\lambda + (n-1)K)|\nabla v|^2+Av^2(A-\lambda)-a\lambda
Av,\nonumber
\end{eqnarray}
where we have used (\ref{d1}) and (\ref{ricci-bound}). Therefore
at $x_0$,
\begin{equation}
0\geq (A-\lambda)|\nabla v|^2+A(A-\lambda)v^2-a\lambda A v
\label{d2}
\end{equation}
and
\[
|\nabla v(x_0)|^2+ \lambda(1+a)v(x_0)^2\leq \frac{a\lambda
v(x_0)}{a\lambda +\epsilon}[\lambda(1+a) +\epsilon]\leq
[\lambda(1+a) +\epsilon].
\]
Finally let $\epsilon\rightarrow 0$. So we have the estimate
(\ref{basic3-n}) in the second case as well.
\end{proof}

We proceed to improve the above bound.

Define a function $Z$ by
\[  Z(t)=\max_{x\in M,t=\sin^{-1} \left(v\left(x\right)/b\right)}
\frac{\left |\nabla v\right |^2}{b^2-v^2}/\lambda.
\]
The estimate in (\ref{basic3-n}) implies that
\begin{equation}                                \label{basic5-n}
Z(t)\leq 1+ a\qquad \textrm{on } [-\sin^{-1}(1/b),
\sin^{-1}(1/b)].
\end{equation}

Throughout this paper, we denote $a/b$ by $c$ and set
\begin{equation}\label{alpha-delta}
\alpha =
\frac12(n-1)K\qquad\textrm{and}\qquad\delta=\alpha/\lambda.
\end{equation}
By (\ref{nK-bound}) we have
\begin{equation}\label{delta-bound}
\delta\leq \frac{n-1}{2n}.
\end{equation}
We have the following conditions on the test function.

\begin{theorem}                                                     \label{thm-barrier-n}
If the function $z:\ [-\sin^{-1} (1/b),\,\sin^{-1} (1/b)]\mapsto
\mathbf{R}^1$ satisfies the following
\begin{enumerate}
 \item $z(t)\geq Z(t) \qquad t\in [-\sin^{-1}(1/b),  \sin^{-1}(1/b)]$,
 \item there exists some $x_0\in M$
       such that at point $t_0=\sin^{-1} (v(x_0)/b)$ \linebreak
       $z(t_0)=Z(t_0)$,
 \item $z(t_0)>0$,
\end{enumerate}
then we have the following
\begin{eqnarray}                    \label{barrier-eq-n}
0&\leq&\frac12z''(t_0)\cos^2t_0 -z'(t_0)\cos t_0\sin t_0 - z(t_0)+1+c\sin t_0 -2\delta \cos^2t_0\nonumber\\
 & & {}-\frac{z'(t_0)}{4z(t_0)}\cos t_0[z'(t_0)\cos t_0 -2z(t_0)\sin t_0 + 2\sin t_0 + 2c].
\end{eqnarray}

\end{theorem}

\begin{corollary}                                                   \label{corollary1}
If in addition to the above conditions 1-3 in Theorem
\ref{thm-barrier-n}, $z'(t_0)\geq 0$ and $1-c\leq z(t_0)\leq 1+a$,
then we have the following
\begin{equation}
0\leq\frac12z''(t_0)\cos^2t_0 -z'(t_0)\cos t_0\sin t_0-z(t_0) +
1+c\sin t_0-2\delta \cos^2t_0.\nonumber
\end{equation}
\end{corollary}

\begin{corollary}                                          \label{corollary2}
If $a=0$, which is defined in (\ref{a-def}), and if  in addition
to the above conditions 1-3 in Theorem \ref{thm-barrier-n},
$z'(t_0)\sin t_0\geq 0$ and $z(t_0)\leq 1$, then we have the
following
\begin{equation}
0\leq\frac12z''(t_0)\cos^2t_0 -z'(t_0)\cos t_0\sin t_0 -z(t_0)+
1-2\delta \cos^2t_0.\nonumber
\end{equation}
\end{corollary}

\begin{proof}[Proof of Theorem \ref{thm-barrier-n}]\quad Define
\[ J(x)=\left\{ \frac{\left |\nabla v\right |^2}{b^2-v^2}
-\lambda z \right\}\cos^2t,
\]
where $t=\sin^{-1}(v(x)/b)$. Then
\[ J(x)\leq 0\quad\textrm{for } x\in M
\qquad \textrm{and} \qquad J(x_0)=0.
\]
If $\nabla v(x_0)=0$ then
\[ 0=J(x_0)=-\lambda z\cos^2 t.
\]
This contradicts the condition 3 in the theorem. Therefore
\[ \nabla v(x_0)\not=0.
\]
Now if $x_0\in M =M\backslash\partial M$ then by the Maximum
Principle, we have
\begin{equation}                                                \label{es1}
\nabla J(x_0)=0\qquad \textrm{and}\qquad \Delta J(x_0)\leq 0.
\end{equation}
If $x_0\in\partial M$, then the weak convexity of $M$, the fact
that $J(x_0)$ is the maximum and an argument in the proof of Lemma
\ref{pre-es-n} imply that $J(x_0)=0$ and $\Delta J(x_0)\leq 0$.
Therefore (\ref{es1}) holds, no matter $x_0\in M
=M\backslash\partial M$ or $x_0\in\partial M$.

$J(x)$ can be rewritten as
\[  J(x)=\frac{1}{b^2}|\nabla v|^2-\lambda z\cos^2t.
\]
Thus (\ref{es1}) is equivalent to
\begin{equation}                                              \label{es2}
\frac{2}{b^2}\sum_{i}v_iv_{ij}\Big|_{x_0}=\lambda\cos t[z' \cos t
-2z\sin t]t_j\Big|_{x_0}
\end{equation}
and
\begin{eqnarray}                                              \label{es3}
0&\geq&\frac{2}{b^2}\sum_{i,j}v_{ij}^2+\frac{2}{b^2}\sum_{i,j}v_iv_{ijj}
 -\lambda (z''|\nabla t|^2+z'\Delta t)\cos^2t \\
 & &+4\lambda z'\cos t\sin t |\nabla t|^2 -
\lambda z\Delta\cos^2t\Big|_{x_0}.\nonumber
\end{eqnarray}
Rotate the local normal frame about $x_0$ such that
$v_1(x_0)\not=0$ and $v_i(x_0)=0$ for $i\geq 2$. Then (\ref{es2})
implies
\begin{equation}                                             \label{es4}
v_{11}\Big|_{x_0}=\frac{\lambda b}{2}(z'\cos t-2z\sin t)
\Big|_{x_0}\quad\text{and}\quad v_{1i} \Big|_{x_0}=0\ \text{for }
i\geq2.
\end{equation}
Now we have
\begin{eqnarray}
|\nabla v|^2
\Big|_{x_0}&=&\lambda b^2z\cos^2t\Big|_{x_0},\nonumber\\
 |\nabla t|^2
\Big|_{x_0}&=&\frac{|\nabla v|^2}{b^2-v^2}=\lambda z
\Big|_{x_0},\nonumber\\
\frac{\Delta v}{b}\Big|_{x_0} &=&\Delta \sin t =\cos t\Delta
t-\sin t |\nabla t|^2
\Big|_{x_0},\nonumber\\
\Delta t\Big|_{x_0}&=&\frac{1}{\cos t}(\sin t|\nabla
t|^2+\frac{\Delta v}{b})
\nonumber\\
 &=&\frac{1}{\cos t} [ \lambda z\sin t-\frac{\lambda}{b}(v+a)] \Big|_{x_0}, \quad\textrm{and}
\nonumber\\
\Delta\cos^2t\Big|_{x_0}&=&\Delta \left(1-\frac{v^2}{b^2}\right)
 =-\frac{2}{b^2}|\nabla v|^2-\frac{2}{b^2}v\Delta v
\nonumber\\
  &=&-2\lambda z\cos^2t+\frac{2}{b^2}\lambda v(v+a)\Big|_{x_0}. \nonumber
\end{eqnarray}
Therefore,
\begin{eqnarray}
& {}&
\frac{2}{b^2}\sum_{i,j}v_{ij}^2\Big|_{x_0}\geq\frac{2}{b^2}v_{11}^2
\nonumber\\
& {}& =\frac{\lambda ^2}{2}(z')^2\cos^2t-2\lambda ^2zz'\cos t\sin
t
      +2\lambda ^2z^2\sin^2 t\Big|_{x_0}\nonumber,
\end{eqnarray}
\begin{eqnarray}
\frac{2}{b^2}\sum_{i,j}v_iv_{ijj}\Big|_{x_0}
&=&\frac{2}{b^2}\left(\nabla v\,\nabla
       (\Delta v)+\textrm{Ric}(\nabla v,\nabla v)\right)\nonumber\\
&\geq& \frac{2}{b^2}(\nabla v\,\nabla (\Delta v)+(n-1)K|\nabla v|^2)\nonumber\\
 &=&-2\lambda^2z\cos^2t+4\alpha \lambda z\cos^2t\Big|_{x_0},\nonumber
\end{eqnarray}
\begin{eqnarray}
&{}&  -\lambda (z''|\nabla t|^2+ z'\Delta t)\cos^2t\Big|_{x_0}\nonumber\\
&{}&=-\lambda^2 zz''\cos^2t-
\lambda^2zz'\cos t\sin t\nonumber\\
&{}&{ }+\frac{1}{b}\lambda^2z'(v +a)\cos t \Big|_{x_0},\nonumber
\end{eqnarray}
and
\begin{eqnarray}
&{}&4\lambda z'\cos t\sin t|\nabla t|^2-\lambda z\Delta
\cos^2t\Big|_{x_0}
\nonumber\\
&{}&=4\lambda^2zz'\cos t\sin
t+2\lambda^2z^2\cos^2t-\frac{2}{b}\lambda^2z\sin t\,(v+a)
\Big|_{x_0}.\nonumber
\end{eqnarray}
Putting these results into (\ref{es3}) we get
\begin{eqnarray}                                                   \label{es5}
0&\geq&-\lambda^2zz''\cos^2t+ \frac{\lambda^2}{2}(z')^2\cos^2t
+\lambda^2z'\cos t\left(z\sin t+c +\sin t\right)
  \nonumber\\
 & & {}+2\lambda^2z^2-2\lambda^2z-2\lambda^2cz\sin t +4\alpha \lambda z\cos^2t
 \Big|_{x_0},
\end{eqnarray}
where we used (\ref{es4}). Now
\begin{equation}                                                    \label{es6}
z(t_0)>0,
\end{equation}
by the condition 3 in the theorem. Dividing two sides of
(\ref{es5}) by $2\lambda^2z\Big|_{x_0}$, we have
\begin{eqnarray}
0&\geq&-\frac12z''(t_0)\cos^2t_0 +\frac12z'(t_0)\cos t_0\left(\sin t_0+\frac{c +\sin t_0}{z(t_0)}\right) +z(t_0) \nonumber\\
 & & {}  -1-c\sin t_0 +2\delta \cos^2t_0\nonumber\\
 & & {}+\frac{1}{4z(t_0)}(z'(t_0))^2\cos^2t_0.\nonumber
\end{eqnarray}
Therefore,
\begin{eqnarray}
0&\geq&-\frac12z''(t_0)\cos^2t_0 + z'(t_0)\cos t_0\sin t_0+z(t_0) -1-c\sin t_0 +2\delta \cos^2t_0\nonumber\\
 & & {}+\frac{z'(t_0)}{4z(t_0)}\cos t_0[z'(t_0)\cos t_0 -2z(t_0)\sin t_0 + 2\sin t_0 + 2c].\nonumber
\end{eqnarray}
\end{proof}

\begin{proof}[Proof of Corollary  \ref{corollary1}] \quad By Condition 2
in the theorem, (\ref{basic5-n}), $|\sin t_0|=|v(t_0)/b|\leq1/b$
and $1-c \leq z(t_0)\leq 1+a$. Thus for $t_0\geq 0$,
\[
-z(t_0)\sin t_0+ \sin t_0 + c  \geq -\sin t_0-a\sin t_0+\sin
t_0+c\geq a(\frac{1}{b}-\sin t_0)\geq 0,
\]
and for $t_0<0$,
\[
-z(t_0)\sin t_0+ \sin t_0 + c  \geq -\sin t_0 + c\sin t_0 +\sin
t_0 +c\geq c(1+\sin t_0)\geq 0.
\]
In any case the last term in the (\ref{barrier-eq-n}) is
non-negative.
\end{proof}

\begin{proof}[Proof of Corollary  \ref{corollary2}]\quad The last term
in the (\ref{barrier-eq-n}) is nonnegative.
\end{proof}

\section{Proof of the Main Result}\label{sec-proofs}
\begin{theorem}                                          \label{thm2-n}
If $a>0$ and $\mu\delta \leq \frac{4}{\pi^2}a$ for a constant
$\mu\in (0,1]$, then under the conditions in Theorem
\ref{main-thm} the first non-zero (closed or Neumann, which
applies) eigenvalue $\lambda$ has the following lower bound
\begin{equation}                                           \label{ling1-1}
\lambda \geq \frac{\pi^2}{d^2} +
\frac{\mu}{2}(n-1)K=\frac{\pi^2}{d^2} + \mu\alpha.
\end{equation}

\end{theorem}

\begin{proof}\quad Let $\mu_{\epsilon}=\mu-\epsilon>0$ for small
positive constant $\epsilon$. Take $b>1$ close to $1$ such that
$\mu_{\epsilon}\delta < \frac{4}{\pi^2}c$. Let
\begin{equation}                                \label{z-def}
z(t)=1+c\eta(t) +\mu_{\epsilon}\delta\xi(t),
\end{equation}
where $\xi$ and $\eta$ are the functions defined by (\ref{xi-def})
and (\ref{eta-def}), respectively.  Let
$\bar{I}=[-\sin^{-1}(1/b),\sin^{-1}(1/b)]$. We claim that
\begin{equation}\label{4.1}
Z(t)\leq z(t)\qquad \textrm{for }t\in \bar{I}.
\end{equation}
By Lemma \ref{xi-lemma} and Lemma \ref{eta-lemma} we have

\begin{eqnarray}
& &{}\frac{1}{2}z''\cos ^2t-z'\cos t\sin t-z
    =-1-c\sin t+ 2\mu_{\epsilon}\delta\cos^2t,          \label{z-eq}\\
& &{}z'(t)> 0\label{z'-geq0}\\
& &{}0<1-\frac{a}{b}=z(-\frac{\pi}{2})\leq z(t)\leq
z(\frac{\pi}{2})=1+\frac{a}{b}\leq 1+a,          \label{z-endpts}
\end{eqnarray}
where (\ref{z'-geq0}) is due to the following.
\begin{eqnarray}
z'(t)=c\eta'(t)+\mu_{\epsilon}\delta\xi'(t)&=&\mu_{\epsilon}\delta\eta'(t)\left(
\frac{c}{\mu_{\epsilon}\delta}+\frac{\xi'(t)}{\eta'(t)}\right)\nonumber\\
&\geq & \mu_{\epsilon}\delta
\eta'(t)(\frac{c}{\mu_{\epsilon}\delta}-\frac{\pi^2}{4}) >
0.\nonumber
\end{eqnarray}
Let $P\in\mathbf{R}^1$ and $t_0\in
[-\sin^{-1}(1/b),\sin^{-1}(1/b)]$ such that
\[ P=\max_{t\in \bar{I}}\left(Z(t)-z(t)\right)=Z(t_0)-z(t_0).
\]
Thus
\begin{equation}\label{4.2}
Z(t)\leq z(t)+P\quad \textrm{for }t\in
\bar{I}\qquad\textrm{and}\qquad Z(t_0)=z(t_0)+P.
\end{equation}
Suppose that $P>0$ Then $z+P$ satisfies the inequality in
Corollary \ref{corollary1} of Theorem \ref{thm-barrier-n}. Then
\begin{eqnarray}
&{}&z(t_0)+P=Z(t_0)\nonumber\\
&{}&\leq  \frac12(z+P)''(t_0)\cos^2 t_0-(z+P)'(t_0)\cos t_0
 \sin t_0+1+c\sin t_0-2\delta \cos^2 t_0\nonumber\\
&{}&=\frac12z''(t_0)\cos^2t_0-z'(t_0)\cos t_0\sin
t_0+1+c\sin t_0-2\delta \cos^2 t_0\nonumber\\
&{}&\leq\frac12z''(t_0)\cos^2t_0-z'(t_0)\cos t_0\sin
t_0+1+c\sin t_0-2\mu_{\epsilon}\delta \cos^2 t_0\nonumber\\
&{}&=z(t_0).\nonumber
\end{eqnarray}
This contradicts the assumption $P>0$. Thus $P\leq 0$ and (\ref{4.1})
must hold. Now we have
\[
|\nabla t|^2\leq\lambda z(t) \qquad \textrm{for } t\in\bar{I},
\]
that is
\begin{equation}\label{4.3}
\sqrt{\lambda}\geq \frac{|\nabla t|}{\sqrt{z(t)}}.
\end{equation}
Let $q_1$ and $q_2$ be two points in $M$ such that $v(q_1)=-1$ and
$v(q_2)=1$ and let $L$ be the minimum geodesic segment between
$q_1$ and $q_2$. We integrate the both sides of (\ref{4.3}) along
$L$ and change variable and let $b\rightarrow 1$. Then
\begin{equation}\label{4.4}
\sqrt{\lambda}d\geq \int_{L}\,\frac{|\nabla
t|}{\sqrt{z(t)}}dl=\int_{-\frac{\pi}{2}}^{\frac{\pi}{2}}
\frac{1}{\sqrt{z(t)}}\,dt \geq \frac{\left(\int_{-\pi/2}^{\pi/2}\
\,dt\right)^\frac32}{(\int_{-\pi/2}^{\pi/2}\ z(t)\,dt)^{\frac12}}
\geq \left( \frac{\pi^3}{\int_{-\pi/2}^{\pi/2}\  z(t)\,dt}
\right)^{\frac12}.
\end{equation}
Square the two sides. Then
\[
\lambda \geq \frac{\pi^3}{d^2\int_{-\pi/2}^{\pi/2} \ z(t)\,dt}.
\]
Now
\[
\int_{-\frac{\pi}{2}}^{\frac{\pi}{2}}\
z(t)\,dt=\int_{-\frac{\pi}{2}}^{\frac{\pi}{2}}\ [1+ a\eta(t)+
\mu_{\epsilon}\delta \xi(t)]\,dt=(1-\mu_{\epsilon}\delta)\pi,
\]
where we used the facts that
$\int_{-\frac{\pi}{2}}^{\frac{\pi}{2}}\ \eta(t)\,dt=0$ since
$\eta$ is an even function, and that $
\int_{-\frac{\pi}{2}}^{\frac{\pi}{2}}\ \xi(t)\,dt=-\pi$ (by
(\ref{xi-int}) in the Lemma \ref{xi-lemma}). Therefore
\[
\lambda \geq
\frac{\pi^2}{(1-\mu_{\epsilon}\delta)d^2}\quad\textrm{and}\quad
\lambda \geq \frac{\pi^2}{d^2} + \mu_{\epsilon}\alpha.
\]
Letting $\epsilon\rightarrow 0$, we get
\[
\lambda \geq \frac{\pi^2}{(1-\mu\delta)d^2}\quad\textrm{and}\quad
\lambda \geq \frac{\pi^2}{d^2} + \mu\alpha.
\]
\end{proof}

\begin{theorem}                                          \label{thm5-n}
If $a=0$, then under the conditions in Theorem \ref{main-thm} the
first non-zero (closed or Neumann, which applies) eigenvalue
$\lambda$ has the following lower bound
\begin{equation}
\lambda \geq \frac{\pi^2}{d^2}+\frac12 (n-1)K.
\end{equation}
\end{theorem}

\begin{proof}\quad Let
\[
y(t)=1+\delta\xi.
\]
By Lemma \ref{xi-lemma},
 for $-\frac{\pi}{2}< t<\frac{\pi}{2}$, we have
\begin{eqnarray}
& &{}\frac{1}{2}y''\cos ^2t-y'\cos t\sin t-y
    =-1+2\delta\cos^2t,\label{21.1}\\
& &{}y'(t)\sin t\geq 0,\qquad \textrm{and} \label{y'-geq0}\\
& &{}y(\pm \frac{\pi}{2})=1 \  \textrm{and }\,
0<y(t)<1.\label{23.3}
\end{eqnarray}
We need only show that $Z(t)\leq y(t)$ on $[-\pi/2,\pi/2]$.  If it
is not true, then there is $t_0$ and a number $P>0$ such that
$P=Z(t_0)-y(t_0)=\max Z(t)-y(t)$. Note that $y(t)+P\geq
1-\frac12(\frac{\pi^2}{4}-1)+P>0$. So $y+P$ satisfies the
inequality in the Corollary \ref{corollary2} in Theorem
\ref{thm-barrier-n}. Therefore
\begin{eqnarray}
&{}&y(t_0)+P=Z(t_0)\nonumber\\
&{}&\leq  \frac12(y+P)''(t_0)\cos^2 t_0-(y+P)'(t_0)\cos t_0
 \sin t_0+1-2\delta \cos^2 t_0\nonumber\\
&{}&=\frac12y''(t_0)\cos^2t_0-y'(t_0)\cos t_0\sin
t_0+1-2\delta \cos^2 t_0\nonumber\\
&{}&=y(t_0).\nonumber
\end{eqnarray}
This contradicts the assumption $P>0$. The rest of the proof is
similar to that of Theorem \ref{thm2-n}, just noticing that
$\delta\leq \frac{n-1}{2n}<\frac12<\frac{4}{\pi^2-4}$.
\end{proof}

\begin{proof}[Proof of Theorem \ref{main-thm} (The Main Result)]\quad
Since $0\leq a<1$, either $a=0$ or $0<a<1$.

If $a=0$, then we apply Theorem \ref{thm5-n} to get the bound with
$\mu=1$,
\[
\lambda \geq \frac{\pi^2}{d^2} + \alpha= \frac{\pi^2}{d^2} +
\frac12(n-1)K.
\]

If $0<a<1$, then there are several cases altogether.

\begin{itemize}
\item(I): $\ $$a\geq\frac{\pi^2}{4}\delta$. \item(II):
$a<\frac{\pi^2}{4}\delta$.
    \begin{itemize}
    \item(II-a): $a\geq 0.765$.
    \item(II-b): $0<a<0.765$.
        \begin{itemize}
        \item(II-b-1): $a\geq 1.53\delta$.
        \item(II-b-2): $a<1.53\delta$.
        \end{itemize}
    \end{itemize}
\end{itemize}

For Case (I): $\ $ $0<a<1$ and $a\geq\frac{\pi^2}{4}\delta$, we
apply Theorem \ref{thm2-n} for $\mu=1$ to get the following lower
bound
\[
\frac{\pi^2}{d^2} + \frac{1}{2}(n-1)K.
\]

For Case (II-a): $0.765\leq a<\frac{\pi^2}{4}\delta$, we apply
Theorem \ref{thm2-n} with $\mu=\frac{4}{\pi^2}\frac{a}{\delta}$
since $(\frac{4}{\pi^2}\frac{a}{\delta})\,\delta\leq
\frac{4}{\pi^2}a$ and $0<\frac{4}{\pi^2}\frac{a}{\delta}<1$. Then
\[
\lambda\geq
\frac{\pi^2}{d^2}+\frac{4}{\pi^2}\frac{a}{\delta}\alpha
=\frac{\pi^2}{d^2}+\frac{4a}{\pi^2}\lambda
\]
Thus
\[
\lambda \geq \frac{1}{1-\frac{4a}{\pi^2}}\,\frac{\pi^2}{d^2}.
\]
On the other hand we have Lichnerowicz-type lower bound
(\ref{delta-bound}),
\[
\lambda \geq \frac{2n}{n-1}\alpha.
\]
The above two estimates give
\[
\lambda \geq \frac{\pi^2}{d^2}
+\frac{4a}{\pi^2}\frac{2n}{n-1}\alpha \geq \frac{\pi^2}{d^2}
+\frac{8(0.765)n}{\pi^2(n-1)}\alpha  \]
\[
> \frac{\pi^2}{d^2}+
\frac{0.62n}{n-1}\alpha >\frac{\pi^2}{d^2}+ \frac{31}{50}\alpha
\]
\[=\frac{\pi^2}{d^2}+ \frac{31}{100}(n-1)K.
\]
The theorem is proved in this case.

For Case (II-b-1): $0<a<0.765$, $a<\frac{\pi^2}{4}\delta$ and
$a\geq 1.53\delta$, we apply Theorem \ref{thm2-n} with with
$\mu=\frac{4}{\pi^2}\frac{a}{\delta}$ since
$(\frac{4}{\pi^2}\frac{a}{\delta})\,\delta\leq \frac{4}{\pi^2}a$
and $0<\frac{4}{\pi^2}\frac{a}{\delta}<1$. Then
\[
\lambda\geq
\frac{\pi^2}{d^2}+\frac{4}{\pi^2}\frac{a}{\delta}\alpha \geq
\frac{\pi^2}{d^2}+\frac{4}{\pi^2}\frac{153}{100}\alpha
\]
\[
>\frac{\pi^2}{d^2}+\frac{31}{50}\alpha
\]
\[
=\frac{\pi^2}{d^2}+\frac{31}{100}(n-1)K,
\]
which is what we want to prove.

For the remaining Case (II-b-2): $0<a<0.765$,
$a<\frac{\pi^2}{4}\delta$ and $a<1.53\delta$, we define a function
$z$ by
\[
z(t)=1+c\eta(t)+(\delta-\sigma c^2)\xi(t) \qquad\textrm{on
}[-\sin^{-1}\frac1b, \sin^{-1}\frac1b],
\]
where
\begin{equation}                                      \label{sigma-def}
\sigma =\frac{\tau}{\left([\,\frac32-\frac{\pi^2}{8}-
(\frac{\pi^2}{32}-\frac16)\frac{153}{100}]\frac{200}{153}-
\frac{(\frac{8}{3\pi}-\frac{\pi}{4})^2}{[-1
+(12-\pi^2)\frac{100}{153}]}\right)c}>0
\end{equation}
and
\begin{equation}                                      \label{tau-def}
\tau =
\frac{2}{3\pi^2}\left(\frac{4}{3(4-\pi)}+\frac{3(4-\pi)}{4}-2\right)>0.
\end{equation}

Let $\bar{I}=[-\sin^{-1}\frac1b, \sin^{-1}\frac1b]$.

 We now show that
\begin{equation}\label{Z-leq-z}
Z(t)\leq z(t)\qquad \textrm{on }\bar{I}.
\end{equation}

If (\ref{Z-leq-z}) is not true, then there exists a constant $P>0$
and $t_0$ such that
\[
Pc^2=\frac{Z(t_0)-z(t_0)}{-\xi(t_0)}=\max_{t\in[-\sin^{-1}\frac1b,
\sin^{-1}\frac1b]}\frac{Z(t)-z(t)}{-\xi(t)}.
\]
Let $w(t)=z(t)-Pc^2\xi(t)=1+c\eta(t) +m\xi(t)$, where $m =
\delta-\sigma c^2 -Pc^2$. Then
\[
Z(t)\leq w(t)\qquad \textrm{on }\bar{I}\qquad \textrm{and}\qquad
Z(t_0)=w(t_0).
\]
By Lemma \ref{z-positive-lemma}, $w(t_0)>0$. So $w$ satisfies
(\ref{barrier-eq-n}) in Theorem \ref{thm-barrier-n},
\[
0\leq -2(\sigma + P)c^2 \cos^2t_0
  -\frac{w'(t_0)}{4w(t_0)}\cos t_0\left(\frac{8c}{\pi}\cos t +4m t\cos
  t\right).
\]
We used (\ref{xi-eq}), (\ref{xi-eq2}), (\ref{eta-eq}) and
(\ref{eta-eq2}) to get the above inequality. Thus
\begin{equation}\label{48.0}
\sigma + P\leq
  -\frac{w'(t_0)}{2c^2w(t_0)}\left(\frac{2c}{\pi} +m t\right)
  =-\frac{\eta'(t_0)}{\pi
  w(t_0)}\left(1+\frac{m\xi'(t_0)}{c\eta'(t_0)}\right)
     \left(1 +\frac{\pi m}{2c}t_0\right).
\end{equation}
The righthand side is not positive for $t_0\geq 0$, by Lemmas
\ref{xi-lemma} and \ref{eta-lemma}. Thus $t_0<0$, and
\begin{eqnarray}
 {}&{}&-\left(1+\frac{m\xi'(t_0)}{c\eta'(t_0)}\right)
     \left(1 +\frac{\pi m}{2c}t_0\right)
\nonumber\\
{}&{}&=\frac{2\xi'(t_0)}{\pi t_0\eta'(t_0)}\left(\frac{\pi
t_0\eta'(t_0)}{2\xi'(t_0)}+\frac{\pi m}{2c}t_0\right)\left(-1
-\frac{\pi m}{2c}t_0\right)\nonumber\\
{}&{}&\leq \frac14\frac{2\xi'(t_0)}{\pi t_0\eta'(t_0)}\left(
\frac{\pi
t_0\eta'(t_0)}{2\xi'(t_0)}-1 \right)^2\nonumber\\
{}&{}&= \frac14\left( \frac{2\xi'(t_0)}{\pi
t_0\eta'(t_0)}+(\frac{2\xi'(t_0)}{\pi
t_0\eta'(t_0)})^{-1}-2\right). \nonumber
\end{eqnarray}
By Lemmas \ref{xi-lemma} and \ref{eta-lemma}, we have
$2(3-\frac{\pi^2}{4})\leq \frac{\xi'(t)}{t}\leq \frac43$ and
$2(\frac{4}{\pi}-1)\leq \eta'(t)\leq \frac{8}{3\pi}$. So
\[
\frac{3(12-\pi^2)}{8}\leq\frac{2\xi'(t_0)}{\pi t_0\eta'(t_0)}\leq
\frac{4}{3(4-\pi)}.
\]
Note that the function $f(t)=t+\frac1t-2$ achieves it maximum on
$[A, B]$ not containing $0$ at an endpoint. Therefore
\[
\left|-\left(1+\frac{m\xi'(t_0)}{c\eta'(t_0)}\right)
     \left(1 +\frac{\pi m}{2c}t_0\right)\right|\leq
     \frac14\left(\frac{4}{3(4-\pi)}+\frac{3}{3(4-\pi)}-2\right).
\]

Now (\ref{48.0}) becomes
\begin{equation}  \label{49}
\sigma + P\leq \frac{\tau}{w(t_0)}.
\end{equation}

On the other hand, by Lemma \ref{z-positive-lemma},
\begin{equation}\label{50}
z(t_0)\geq
\left([\,\frac32-\frac{\pi^2}{8}-(\frac{\pi^2}{32}-\frac16)\frac{153}{100}]\frac{200}{153}-
\frac{(\frac{8}{3\pi}-\frac{\pi}{4})^2}{[-1
+(12-\pi^2)\frac{100}{153}]}\right)c=\frac{\tau}{\sigma}>0.
\end{equation}
Since $-P\xi(t_0)\geq 0$, we have $w(t_0)\geq z(t_0)$. This fact,
(\ref{49}) and (\ref{50}) imply that for $P>0$
\[
\sigma + P<\sigma,
\]
which is impossible.

Therefore we have the estimate (\ref{Z-leq-z}). Now we proceed as
in the proof of Theorem \ref{thm2-n}. We get the following
\[
\lambda d^2\geq \frac{\pi^3}{\pi[1-(\delta-\sigma c^2)]}.
\]
Since $\delta -\sigma c^2>0.625\delta$ by Lemma
\ref{z-positive-lemma}, we have
\[
\lambda \geq \frac{1}{[1-(\delta-\sigma
c^2)]}\frac{\pi^2}{d^2}\geq
\frac{1}{[1-0.625\delta]}\frac{\pi^2}{d^2}
\]
and
\[
\lambda \geq
\frac{\pi^2}{d^2}+0.625\alpha>\frac{\pi^2}{d^2}+\frac{31}{100}(n-1)K.
\]

If $n=2$ then we can get even better result.

If $a=0$, then we apply Theorem \ref{thm5-n} to get the lower
bound $\frac{\pi^2}{d^2}+\frac12(n-1)K$.

If $a\geq\frac{\pi^2}{4}\delta$, then we apply Theorem
\ref{thm2-n} to get the lower bound
$\frac{\pi^2}{d^2}+\frac12(n-1)K$.

If $a<\frac{\pi^2}{4}\delta$ and $n=2$, then $a$ satisfies
\[a\leq
\frac{(12-\pi^2)n+\pi^2-4}{8n}.
\]
Otherwise that
$a<\frac{\pi^2}{4}\delta$, $a> \frac{(12-\pi^2)n+\pi^2-4}{8n}$ and
$\delta\leq \frac{n-1}{2n}$ would yield
\[\frac{(12-\pi^2)n+\pi^2-4}{8n}<a< \frac{\pi^2}{4}\delta\leq
\frac{\pi^2}{4}\frac{n-1}{2n}=\frac{\pi^2}{16}=\frac{\pi^2(n-1)}{8n}.\]
We do know the following opposite inequality holds for $n=2$,
\[
\frac{(12-\pi^2)n+\pi^2-4}{8n}=\frac{20-\pi^2}{16}>\frac{\pi^2}{16}=\frac{\pi^2(n-1)}{8n}.
\]
Therefore we may apply Theorem \ref{thm3-n} to get the the lower
bound stated in the theorem, which is the least of the three lower
bounds.
\end{proof}

We now present a Lemma that is used in the proof of the Theorem
\ref{main-thm}.

\begin{lemma}\label{z-positive-lemma}
If $a<1.53\delta$ and $0<a<0.765$ then
\[
z(t)=1+c\eta(t)+\delta\xi(t)
\]
\[
\geq
\left([\,\frac32-\frac{\pi^2}{8}-(\frac{\pi^2}{32}-\frac16)\frac{153}{100}]\frac{200}{153}-
\frac{(\frac{8}{3\pi}-\frac{\pi}{4})^2}{[-1
+(12-\pi^2)\frac{100}{153}]}\right)c>0,
\]
for $t\in[-\pi/2, \pi/2]$ and
\[
\delta-\sigma c^2 \approx 0.625162283437>0.625\delta,
\]
where $c=a/b$ and $b>1$ is any constant and $\sigma$ is the
constant in (\ref{sigma-def}).
\end{lemma}
\begin{proof}\quad By Lemmas \ref{r-lemma}, Lemma \ref{xi-lemma} and \ref{eta-lemma},
the function $z$ on $[-\pi/2, \pi/2]$ has a unique critical point
$t_1\in (-\pi/2, 0)$  if $0<a<\frac{\pi^2}{4}\delta$ and $z$ is
decreasing on $[-\pi/2, t_1]$ and increasing on $[t_1, \pi/2]$.
Therefore
\[
\min_{[-\pi/2, \pi/2]}z=\min_{[-\pi/2, 0]}z=z(t_1).
\]
So we need only consider the restricted function $z|_{[-\pi/2,
0]}$ for the minimum.

Now first consider the Taylor expansion of $\xi$ at $0$ for
$t\in[-\pi/2, 0]$. By Lemma \ref{xi-lemma},
$\xi(0)=-\frac{\pi^2}{4}+1$, $\xi'(0)=0$ and
$\xi''(0)=2(3-\frac{\pi^2}{4})$ and $\xi'''(t)<0$ on $(-\pi/2,
0)$.

Thus

\[
\xi(t)=\xi(0)
+\xi'(0)+\frac{\xi''(0)}{2!}t^2+\frac{\xi'''(t_2)}{2!}t^3
\]
\[
\geq\xi(0) +\xi'(0)+\frac{\xi''(0)}{2!}t^2
\]
\[
=-(\frac{\pi^2}{4}-1)+(3-\frac{\pi^2}{4})t^2,
\]
where $t_2$ is a constant in $(t, 0)$. Similarly, using the data
$\eta(-\pi/2)=-1$, $\eta'(-\pi/2)=\frac{8}{3\pi}$ and
$\eta'''(t)>0$ on $(-\pi/2, 0)$ (actually on $[-\pi/2, \pi/2]$),
and the Taylor expansion of $\eta$ at $-\pi/2$,  we have for $t\in
[-\pi/2, 0]$,
\[
\eta(t)=\eta(-\frac{\pi}{2})+\eta'(-\frac{\pi}{2})(t+\frac{\pi}{2})
+\frac{\eta''(-\frac{\pi}{2})}{2!}(t+\frac{\pi}{2})^2
+\frac{\eta''(t_3)}{3!}(t+\frac{\pi}{2})^3
\]
\[
\geq \eta(-\frac{\pi}{2})+\eta'(-\frac{\pi}{2})(t+\frac{\pi}{2})
+\frac{\eta''(-\frac{\pi}{2})}{2!}(t+\frac{\pi}{2})^2
\]
\[
=-1+\frac{8}{3\pi}(t+\frac{\pi}{2})-\frac{1}{4}(t+\frac{\pi}{2})^2
\]
\[=-(\frac{\pi^2}{16}-\frac13)+(\frac{8}{3\pi}-\frac{\pi}{4})t-\frac14
t^2,
\]
where $t_3$ is some constant in $(-\pi/2,t)$. Therefore on
$[-\pi/2,0]$,
\[
z(t)=1+c\eta(t)+\delta\xi(t)
\]
\[
\geq
1-(\frac{\pi^2}{16}-\frac13)c-(\frac{\pi^2}{4}-1)\delta+(\frac{8}{3\pi}-\frac{\pi}{4})ct+[-\frac14c
+(3-\frac{\pi^2}{4})\delta] t^2.
\]

Let $\nu=1.53$ and $a_0=0.765$. That $a\leq \nu\delta$ implies
$c=a/b<\nu \delta$, where $b>1$ is a constant.  Using conditions
(\ref{delta-bound}) $\delta\leq \frac{n-1}{2n}<\frac12$ and $a\leq
a_0$, we get

\[
1-(\frac{\pi^2}{16}-\frac13)c-(\frac{\pi^2}{4}-1)\delta
\]
\[
\geq1-(\frac{\pi^2}{16}-\frac13)\nu\delta-(\frac{\pi^2}{4}-1)\delta
\]
\[
\geq \frac32-\frac{\pi^2}{8}-(\frac{\pi^2}{32}-\frac16)\nu
\]
\[
>
\left(\frac32-\frac{\pi^2}{8}-(\frac{\pi^2}{32}-\frac16)\nu\right)\frac{1}{a_0}c
\]

and

\[
1+c\eta(t)+\delta\xi(t)
\]
\[
\geq
\left(\frac32-\frac{\pi^2}{8}-(\frac{\pi^2}{32}-\frac16)\nu\right)\frac{1}{a_0}c
+(\frac{8}{3\pi}-\frac{\pi}{4})ct+[-\frac14c
+(3-\frac{\pi^2}{4})\frac{1}{\nu}c] t^2
\]
\[
=\left([\,\frac32-\frac{\pi^2}{8}-(\frac{\pi^2}{32}-\frac16)\nu]\frac{1}{a_0}
+(\frac{8}{3\pi}-\frac{\pi}{4})t+[-\frac14
+(3-\frac{\pi^2}{4})\frac{1}{\nu}] t^2\right)c
\]
\[
\geq
\left([\,\frac32-\frac{\pi^2}{8}-(\frac{\pi^2}{32}-\frac16)\nu]\frac{1}{a_0}-
\frac{(\frac{8}{3\pi}-\frac{\pi}{4})^2}{4[-\frac14
+(3-\frac{\pi^2}{4})\frac{1}{\nu}]}\right)c
\]
\[
\geq
\left([\,\frac32-\frac{\pi^2}{8}-(\frac{\pi^2}{32}-\frac16)\nu]\frac{1}{a_0}-
\frac{(\frac{8}{3\pi}-\frac{\pi}{4})^2}{[-1
+(12-\pi^2)\frac{1}{\nu}]}\right)c>0.5433>0.
\]
Let $\tau$ be the constant in (\ref{tau-def}). Then
\[
\sigma c^2 =\frac{\tau
c}{\left([\,\frac32-\frac{\pi^2}{8}-(\frac{\pi^2}{32}-\frac16)\nu]\frac{1}{a_0}-
\frac{(\frac{8}{3\pi}-\frac{\pi}{4})^2}{[-1
+(12-\pi^2)\frac{1}{\nu}]}\right)c},
\]
\[
\leq\frac{\tau\nu
\delta}{\left([\,\frac32-\frac{\pi^2}{8}-(\frac{\pi^2}{32}-\frac16)\nu]\frac{1}{a_0}-
\frac{(\frac{8}{3\pi}-\frac{\pi}{4})^2}{[-1
+(12-\pi^2)\frac{1}{\nu}]}\right)c}\approx 0.374837516563\delta
\]
and
\[
\delta-\sigma c^2 >0.625\delta.
\]
\end{proof}

\begin{theorem}                                         \label{thm3-n}
If $0<a<\frac{\pi^2}{4}\delta$ and $a\leq
\frac{(12-\pi^2)n+\pi^2-4}{8n}$, then under the conditions in
Theorem \ref{main-thm} the first non-zero (closed or Neumann,
which applies) eigenvalue has the following lower bound
\[
\lambda \geq \frac{\pi^2}{d^2} + \frac{\mu}{2}(n-1)K,
\]
where
\begin{equation}       \label{mu-p}
\mu=1-\sqrt{\frac{\pi^2}{6(\pi^2-4)}\left(\frac{4}{3(4-\pi)}
+\frac{3(4-\pi)}{4}-2\right)}\approx (0.765\cdots)>3/4.
\end{equation}
\end{theorem}

\begin{proof}\quad
The proof is similar to that of Case (II)-b-2 in the proof of
Theorem \ref{main-thm}. Clearly, we have
$c<\frac{\pi^2}{4}\delta$, where $c=a/b$ with constant $b>1$. Let
\[
z=1+c\eta +(\delta-\tilde{\sigma} c^2)\xi  \qquad\textrm{on
}[-\sin^{-1}\frac1b, \sin^{-1}\frac1b],
\]
where $\xi$ and
$\eta$ are functions defined in (\ref{xi-def}) and (\ref{eta-def})
respectively, $\tau$ is the constant in (\ref{tau-def}) and
\begin{equation}                                 \label{sigma-tilde-def}
\tilde{\sigma} = \frac{-[1-c-(\frac{\pi^2}{4}-1)\delta
]+\sqrt{[1-c-(\frac{\pi^2}{4}-1)\delta
]^2+4(\frac{\pi^2}{4}-1)\tau c^2}}{2(\frac{\pi^2}{4}-1)c^2}.
\end{equation}

We prove that
\[
Z(t)\leq z(t)\qquad \textrm{on }[-\sin^{-1}\frac1b,
\sin^{-1}\frac1b].
\]
If it is not true, then there exists a constant $P>0$ and $t_0$
such that
\[
Pc^2=\frac{Z(t_0)-z(t_0)}{-\xi(t_0)}=\max_{t\in[-\sin^{-1}\frac1b,
\sin^{-1}\frac1b]}\frac{Z(t)-z(t)}{-\xi(t)}.
\]
Let $\bar{I}=[-\sin^{-1}\frac1b, \sin^{-1}\frac1b]$ and
$w(t)=z(t)-Pc^2\xi(t)=1+c\eta(t) +m\xi(t)$, where $m =
\delta-\tilde{\sigma} c^2 -Pc^2$.   Then
\begin{equation}\label{Z-z-0}
Z(t)\leq w(t) \qquad \textrm{on }\bar{I}\qquad \textrm{and}\qquad
Z(t_0)=w(t_0).
\end{equation}

We want to show that $w(t_0)>0$. In order to do that, we now show
that $m>0$ first.
\begin{lemma}                           \label{z0-lemma}
$Z(t)\leq 1+ c\eta(t)\quad\textrm{on
}[-\sin^{-1}\frac1b,\sin^{-1}\frac1b]$.
\end{lemma}
\begin{proof}[Proof of Lemma \ref{z0-lemma}]\quad If it is not true,
then there exist $t_0$ and constant $P$ such that
$P=Z(t_0)-[1+c\eta (t_0)]= \max \left(Z(t)-[1+c\eta(t)]\right)$.
Thus $1+c\eta +P$ satisfies the inequality in Corollary 1 of the
Theorem \ref{thm-barrier-n}. Therefore
\begin{eqnarray}
&{}&1+\eta(t_0)+P=Z(t_0)\nonumber\\
&{}&\leq  \frac12(1+\eta +P)''(t_0)\cos^2 t_0-(1+\eta+P)'(t_0)\cos
t_0
 \sin t_0+1+c\sin t_0-2\delta \cos^2 t_0\nonumber\\
&{}&=\frac12\eta''(t_0)\cos^2t_0-\eta'(t_0)\cos t_0\sin
t_0+1+c\sin t_0-2\delta \cos^2 t_0\nonumber\\
&{}&=1+\eta(t_0)-2\delta\cos^2t_0\nonumber\\
&{}&\leq 1+\eta(t_0).\nonumber
\end{eqnarray}
This contradicts the assumption $P>0$. The proof of the lemma is
completed.
\end{proof}

Lemma \ref{z0-lemma} implies that $w(t_0)=1+c\eta(t_0)+m
\xi(t_0)=Z(t_0)\leq 1+c\eta(t_0)$. Thus $m \xi(t_0)\leq 0$ and
$m=\delta -\tilde{\sigma} c^2 -Pc^2\geq 0$.

We now show that $w(t_0)>0$. By the fact $m\geq 0$, Lemmas
\ref{xi-lemma} and \ref{eta-lemma}, we have
\begin{eqnarray}
w(t)&\geq& 1-c-(\frac{\pi^2}{4}-1)( \delta- \tilde{\sigma}
c^2-Pc^2)\nonumber\\
{}&>&1-c-(\frac{\pi^2}{4}-1)\delta
+(\frac{\pi^2}{4}-1)(\tilde{\sigma}
c^2+Pc^2)\nonumber\\
{}&>&1-c-(\frac{\pi^2}{4}-1)\delta+(\frac{\pi^2}{4}-1)\tilde{\sigma}
c^2>1-c-(\frac{\pi^2}{4}-1)\delta.\label{48.2-0}
\end{eqnarray}
We claim that if $a\leq \frac{(12-\pi^2)n+\pi^2-4}{8n}$ then
\[w(t)>1-c-(\frac{\pi^2}{4}-1)\delta>0.\]

In fact, (\ref{48.2-0}), (\ref{delta-bound}), and $a\leq
\frac{(12-\pi^2)n+\pi^2-4}{8n}$ imply that
\[
w(t)>1-c-(\frac{\pi^2}{4}-1)\delta>1-a-(\frac{\pi^2}{4}-1)\delta
\]
\[
\geq
1-\frac{(12-\pi^2)n+\pi^2-4}{8n}-(\frac{\pi^2}{4}-1)\frac{n-1}{2n}=0.
\]

Therefore $w(t)>0$ and $\tilde{\sigma}>0$. Now (\ref{Z-z-0}) and
the fact $w(t_0)>0$ imply that $w$ satisfies (\ref{barrier-eq-n})
in Theorem \ref{thm-barrier-n}. So
\[
0\leq -2(\tilde{\sigma} + P)c^2 \cos^2t_0
  -\frac{w'(t_0)}{4w(t_0)}\cos t_0\left(\frac{8c}{\pi}\cos t +4m t\cos t\right),\nonumber
\]
where we used (\ref{xi-eq}), (\ref{xi-eq2}), (\ref{eta-eq}) and
(\ref{eta-eq2}) to get the above inequality. Thus
\begin{equation}\label{48.0-0}
\tilde{\sigma} + P\leq
  -\frac{w'(t_0)}{2c^2w(t_0)}\left(\frac{2c}{\pi} +m t\right)
  =-\frac{\eta'(t_0)}{\pi w(t_0)}\left(1+\frac{m\xi'(t_0)}{c\eta'(t_0)}\right)
     \left(1 +\frac{\pi m}{2c}t_0\right).
\end{equation}
The righthand side is not positive as $t_0\geq 0$, by Lemmas
\ref{xi-lemma} and \ref{eta-lemma}. Thus $t_0<0$. It is showed in
the proof of Case (II)-b-2 of the proof of Theorem \ref{main-thm}
that
\[
\left|-\left(1+\frac{m\xi'(t_0)}{c\eta'(t_0)}\right)
     \left(1 +\frac{\pi m}{2c}t_0\right)\right|\leq
     \frac14\left(\frac{4}{3(4-\pi)}+\frac{3}{3(4-\pi)}-2\right).
\]

Therefore (\ref{48.0-0}) becomes
\begin{equation}  \label{49-0}
\tilde{\sigma} + P\leq \frac{\tau}{w(t_0)}.
\end{equation}

Now taking (\ref{48.2-0}) into (\ref{49-0}), we get
\[
\tilde{\sigma}  + P \leq
\frac{\tau}{w(t_0)}\leq\frac{\tau}{1-c-(\frac{\pi^2}{4}-1)\delta +
(\frac{\pi^2}{4}-1)\tilde{\sigma} c^2}=\tilde{\sigma}.
\]
This contradicts $P> 0$. The last equality is due to the fact that
$\tilde{\sigma}$ is the positive solution of the quadratic
equation
\[
-\tau+(1-c)\tilde{\sigma}-(\frac{\pi^2}{4}-1)\delta\tilde{\sigma}
+ (\frac{\pi^2}{4}-1)\tilde{\sigma}^2 c^2=0.
\]
Therefore
\[ Z(t)\leq z(t) = 1+ c\eta(t)+(\delta-\tilde{\sigma}
c^2)\xi.
\]

Note that $\sqrt{A+B}\leq \sqrt{A} +\sqrt{B}$. By the conditions
$c< \frac{\pi^2}{4}\delta$ and $1-c-(\frac{\pi^2}{4}-1)\delta>0$,
we have
\[
\tilde{\sigma} c^2=\frac{-[1-c-(\frac{\pi^2}{4}-1)\delta ]+
\sqrt{[1-c-(\frac{\pi^2}{4}-1)\delta ]^2 +4(\frac{\pi^2}{4}-1)\tau
c^2}}{2(\frac{\pi^2}{4}-1)}
\]
\[
\leq  \frac{c\sqrt{(\frac{\pi^2}{4}-1)\tau }}{(\frac{\pi^2}{4}-1)}
\leq\frac{\frac{\pi^2}{4}\delta
\sqrt{(\frac{\pi^2}{4}-1)\tau}}{(\frac{\pi^2}{4}-1)}
=\frac{\pi^2\delta}{2} \sqrt{\frac{\tau}{\pi^2-4}}\approx
(0.235\cdots)\delta,
\]
and
\[
\delta-\tilde{\sigma} c^2\geq \left(1-\frac{\pi^2}{2}
\sqrt{\frac{\tau}{\pi^2-4}}\right)\delta=\mu\delta,
\]
where $\mu$ is the constant in (\ref{mu-p}). Proceeding further as
in the proof of Theorem \ref{thm2-n}, we get
\[
\lambda\geq
\frac{1}{1-(\delta-\tilde{\sigma}c^2)}\frac{\pi^2}{d^2}\geq
\frac{1}{1-\mu\delta}\frac{\pi^2}{d^2}
\]
and
\[
\lambda\geq \frac{\pi^2}{d^2}+\mu\alpha.
\]
\end{proof}

\section{Functions}\label{sec-functions}

We study the functions that are used for the construction of the
test functions.

\begin{lemma}                           \label{xi-lemma}
Let
\begin{equation}                        \label{xi-def}
\xi(t)=\frac{\cos^2t+2t\sin t\cos t +t^2-\frac{\pi^2}{4}}{\cos^2t}
\qquad \textrm{on}\quad [-\frac{\pi}{2},\frac{\pi}{2}\,].
\end{equation}
Then the function $\xi$  satisfies the following
\begin{eqnarray}
& &{}\frac{1}{2}\xi''\cos ^2t-\xi'\cos t\sin t-\xi
    =2\cos^2t\quad \textrm{in }(-\frac{\pi}{2},\frac{\pi}{2}\,),          \label{xi-eq}\\
& &{}\xi'\cos t -2\xi\sin t =4t\cos t  \quad \textrm{in }(-\frac{\pi}{2},\frac{\pi}{2}\,),                     \label{xi-eq2}\\
& &{}\int_0^{\frac{\pi}{2}}\xi(t)\, dt= -\frac{\pi}{2} ,         \label{xi-int}\\
& &{}1-\frac{\pi^2}{4}=\xi(0)\leq\xi(t)\leq\xi(\pm
\frac{\pi}{2})=0\quad
\textrm{on }[-\frac{\pi}{2},\frac{\pi}{2}\,],                                \nonumber\\
& &{} \xi' \textrm{ is increasing on }
[-\frac{\pi}{2},\frac{\pi}{2}\,] \textrm{ and }
\xi'(\pm \frac{\pi}{2}) =\pm \frac{2\pi}{3},                     \nonumber\\
& &{}\xi'(t)< 0 \textrm{ on }(-\frac{\pi}{2},0)\textrm{ and \ }
\xi'(t)>0 \textrm{ on }(0,\frac{\pi}{2}\,), \nonumber \\
& &{}\xi''(\pm\frac{\pi}{2})=2, \ \xi''(0)=2(3-\frac{\pi^2}{4})
\textrm{ and \ } \xi''(t)> 0 \textrm{ on }
[-\frac{\pi}{2},\frac{\pi}{2}\,], \nonumber\\
& &{}(\frac{\xi'(t)}{t})'>0 \textrm{ on } (0,\pi/2\,)\textrm{ and
\ } 2(3-\frac{\pi^2}{4})\leq \frac{\xi'(t)}{t}\leq \frac43
\textrm{ on } [-\frac{\pi}{2},\frac{\pi}{2}\,],\nonumber \\
& &{}\xi'''(\frac{\pi}{2})=\frac{8\pi}{15}, \xi'''(t)< 0 \textrm{
on }(-\frac{\pi}{2},0) \textrm{ and \ }  \xi'''(t)>0 \textrm{ on
}(0,\frac{\pi}{2}\,). \nonumber
\end{eqnarray}
\end{lemma}
\begin{proof}\quad For convenience, let $q(t)= \xi'(t)$, i.
e.,
\begin{equation}                                       \label{q-def}
q(t) = \xi'(t) = \frac{2(2t\cos t +t^2\sin t +\cos^2 t \sin t
-\frac{\pi^2}{4}\sin t)}{\cos^3 t}.
\end{equation}
Equation (\ref{xi-eq}) and the values $\xi(\pm \frac{\pi}{2})=0$,
$\xi(0)=1-\frac{\pi^2}{4}$ and $\xi'(\pm \frac{\pi}{2}) =\pm
\frac{2\pi}{3}$ can be verified directly from (\ref{xi-def}) and
(\ref{q-def}) .  The values of $\xi''$ at $0$ and $\pm
\frac{\pi}{2}$ can be computed via (\ref{xi-eq}). By
(\ref{xi-eq2}), $(\xi(t)\cos^2 t)' =4t\cos^2 t$. Therefore
\newline $\xi(t)\cos^2 t=\int_{\frac{\pi}{2}}^t \ 4s\cos^2 s\,ds$,
and
\[
\int_{-\frac{\pi}{2}}^{\frac{\pi}{2}}\
\xi(t)\,dt=2\int_0^{\frac{\pi}{2}}\
\xi(t)\,dt=-8\int_0^{\frac{\pi}{2}}\left( \frac{1}{\cos^2(t)}
\int_t^{\frac{\pi}{2}}\ s\cos^2s\,ds\right)\,dt
\]
\[
=-8\int_0^{\frac{\pi}{2}}\left(\int_0^s\
\frac{1}{\cos^2(t)}\,dt\right)\ s\cos^2s\,ds
=-8\int_0^{\frac{\pi}{2}}\ s\cos s\sin s\,ds=-\pi.
\]
It is easy to see that $q$ and $q'$ satisfy the following
equations
\begin{equation}                                         \label{q-eq}
\frac12 q''\cos t -2q'\sin t -2q\cos t = -4 \sin t,
\end{equation}
and
\begin{equation}                                       \label{q'-eq}
\frac{\cos^2 t}{2(1+\cos^2 t)}(q')''-\frac{2\cos t\sin t}{1+\cos^2
t}(q')'-2(q')=-\frac{4}{1+\cos^2 t}.
\end{equation}
The last equation implies $q'=\xi''$ cannot achieve its
non-positive local minimum at a point in $(-\frac{\pi}{2},
\frac{\pi}{2})$. On the other hand, $\xi''(\pm\frac{\pi}{2})=2$,
by equation (\ref{xi-eq}), $\xi(\pm \frac{\pi}{2})=0$ and
$\xi'(\pm \frac{\pi}{2})=\pm \frac{2\pi}{3}$. Therefore
$\xi''(t)>0$ on $[-\frac{\pi}{2},\frac{\pi}{2}]$ and $\xi'$ is
increasing. Since $\xi'(t)=0$, we have $\xi'(t)< 0$ on
$(-\frac{\pi}{2},0)$ and $\xi'(t)>0$ on $(0,\frac{\pi}{2})$.
Similarly, from the equation
\begin{eqnarray}                                           \label{q''-eq}
&\frac{\cos^2 t}{2(1+\cos^2 t)}(q'')'' -\frac{\cos t\sin t
(3+2\cos^2 t)}{(1+\cos^2 t)^2}(q'')' -\frac{2(5\cos^2 t+\cos^4 t)}{(1+\cos^2 t)^2}(q'') \nonumber\\
&=-\frac{8\cos t\sin t}{(1+\cos^2 t)^2}
\end{eqnarray}
we get the results in the last line of the lemma.

Set $h(t)=\xi''(t)t-\xi'(t)$. Then $h(0)=0$  and $h'(t)=
\xi'''(t)t>0$ in $(0,\frac{\pi}{2})$. Therefore
$(\frac{\xi'(t)}{t})'=\frac{h(t)}{t^2}>0$ in $(0,\frac{\pi}{2})$.
Note that $\frac{\xi'(-t)}{-t}= \frac{\xi'(t)}{t}$,
$\frac{\xi'(t)}{t}|_{t=0}=\xi''(0)=2(3-\frac{\pi^2}{4})$ and
$\frac{\xi'(t)}{t}|_{t=\pi/2}=\frac43$. This completes the proof
of the lemma.
\end{proof}

\begin{lemma}                                       \label{eta-lemma}
Let
\begin{equation}                                    \label{eta-def}
 \eta(t)=\frac{\frac{4}{\pi}t+\frac{4}{\pi}\cos t\sin t-2\sin t}{\cos^2t}
 \qquad
\textrm{on}\quad [-\frac{\pi}{2},\frac{\pi}{2}\,].
\end{equation}
Then the function $\eta$ satisfies the following
\begin{eqnarray}
& &{}\frac{1}{2}\eta''\cos ^2t-\eta'\cos t\sin t-\eta
     =-\sin t\qquad \textrm{in \  }(-\frac{\pi}{2},\frac{\pi}{2}\,),              \label{eta-eq}\\
& &{}\eta'\cos t -2\eta \sin t =\frac{8}{\pi}\cos t -2  \qquad \textrm{in \  }(-\frac{\pi}{2},\frac{\pi}{2}\,),                         \label{eta-eq2}\\
& &{}-1=\eta(-\frac{\pi}{2})\leq\eta(t)\leq\eta(\frac{\pi}{2})=1 \qquad \textrm{on \ }[-\frac{\pi}{2},\frac{\pi}{2}\,],\nonumber\\
& &{}0<2(\frac{4}{\pi}-1)=\eta'(0)\leq \eta'(t)\leq
\eta'(\pm \frac{\pi}{2}) =\frac{8}{3\pi} \qquad \textrm{on \ }[-\frac{\pi}{2},\frac{\pi}{2}\,],\nonumber\\
& &{}-1/2=\eta''(-\frac{\pi}{2})\leq \eta''(t)\leq
\eta''(\frac{\pi}{2}) =1/2\qquad \textrm{on \ }[-\frac{\pi}{2},\frac{\pi}{2}\,],\nonumber\\
& &{}\eta'''(t)>0 \textrm{ \ on \ }[-\frac{\pi}{2},
\frac{\pi}{2}]\quad \textrm{and}\quad \eta'''(\pm
\frac{\pi}{2})=\frac{32}{15\pi}.\nonumber
\end{eqnarray}
\end{lemma}
\begin{proof}\quad Let $p(t)= \eta'(t)$, i.e.,
\begin{equation}                                                \label{p-def}
p(t) = \eta'(t) = \frac{2(\frac{4}{\pi}\cos t+\frac{4}{\pi}t\sin
t-\sin^2t-1)} {\cos^3t}.
\end{equation}
Equation (\ref{eta-eq}), $\eta(\pm \frac{\pi}{2})=\pm 1$,
$\eta'(0)=2(\frac{4}{\pi}-1)$ and $\eta'(\pm \frac{\pi}{2})
=\frac{8}{3\pi}$ can be verified directly.  We get $\eta''(\pm
\frac{\pi}{2}) =\pm 1/2$ from the above values and equation
(\ref{eta-eq}). By (\ref{eta-eq}), $q=\eta'$, $q'=\eta''$ and
$p''=\eta'''$ satisfy the following equations in
$(-\frac{\pi}{2},\frac{\pi}{2})$
\begin{equation}                                        \label{p-eq}
\frac12 p''\cos t -2p'\sin t -2p\cos t = -1,
\end{equation}
\[
\frac{\cos^2 t}{2(1+\cos^2 t)}p'''-\frac{2\cos t\sin t}{1+\cos^2
t}p''-2p'=-\frac{\sin t}{1+\cos^2 t},
\]
and
\begin{eqnarray}                                           \label{p''-eq}
&\frac{\cos^2 t}{2(1+\cos^2 t)}(p'')'' -\frac{\cos t\sin t
(3+2\cos^2 t)}{(1+\cos^2 t)^2}(p'')' -\frac{2(5\cos^2 t+\cos^4 t)}{(1+\cos^2 t)^2}(p'') \nonumber\\
&=-\frac{\cos t(2+\sin t)}{(1+\cos^2 t)^2}.
\end{eqnarray}
The coefficient of $(p'')$ in (\ref{p''-eq}) is obviously negative
in $(-\frac{\pi}{2}, \frac{\pi}{2})$ and the righthand side of
(\ref{p''-eq}) is also negative. So $p''$ cannot achieve its
non-positive local minimum at a point in $(-\frac{\pi}{2},
\frac{\pi}{2})$. On the other hand,
$p''(\frac{\pi}{2})=\frac{32}{15\pi}>0$ (see the proof below),
$p''(t)
> 0$ on $[-\frac{\pi}{2}, \frac{\pi}{2}]$. Therefore $p'$ is increasing
and $-1/2=p'(-\frac{\pi}{2})\leq p'(t)\leq p'(\frac{\pi}{2})
=1/2$. Note that $p'(0)= 0$ ($p'$ is an odd function). So $p'(t)>
0$ on $(0,\frac{\pi}{2})$ and $p$ is increasing on
$[0,\frac{\pi}{2}\,]$. Therefore $2(4/\pi-1)=p(0)\leq
p(t)=\eta'(t)\leq p(\frac{\pi}{2}) =\frac{8}{3\pi}$ on
$[0,\frac{\pi}{2}]$, and on $[-\frac{\pi}{2}, \frac{\pi}{2}]$
since $p$ is an even function.
We now show that $p(\frac{\pi}{2})=\frac{8}{3\pi}$,
$p'(\frac{\pi}{2})=1/2$ and $p''(\frac{\pi}{2})=\frac{32}{15\pi}$.
The first is from a direct computation by using (\ref{p-def}). By
(\ref{eta-eq}),
\[
\frac12p'(\frac{\pi}{2})=\frac12\eta''(\frac{\pi}{2}) =
\lim_{t\rightarrow \frac{\pi}{2}^-}\frac{\eta'(t)\cos t\sin t +
\eta(t)-\sin t}{\cos^2 t}=-\frac12[\eta''(\frac{\pi}{2})-1].
\]
So $p'(\frac{\pi}{2})=1/2$. Similarly, by (\ref{p-eq}),
\[
\frac12p''(\frac{\pi}{2}) =\lim_{t\rightarrow
\frac{\pi}{2}-}\frac{2p'(t)\sin t-1}{\cos t}  +
2p(\frac{\pi}{2})=-2p''(\frac{\pi}{2}) + \frac{16}{3\pi}
\]
Thus $p''(\frac{\pi}{2})=\frac{32}{15\pi}$.
\end{proof}

\begin{lemma}                                                   \label{r-lemma}
The function \ $r(t)=\xi'(t)/\eta'(t)$ is an increasing function
on $[-\frac{\pi}{2}, \frac{\pi}{2}]$, i.e., $r'(t)>0$, and $|r(t)|
\leq \frac{\pi^2}{4}$ holds on $[-\frac{\pi}{2}, \frac{\pi}{2}]$.
\end{lemma}
\begin{proof}\quad Let $p(t) =\eta'(t)$ as in (\ref{p-def})
and $q(t)=\xi'(t)$. Then $r(t)=q(t)/p(t)$. It is easy to  verify
that $ r(\pm \frac{\pi}{2})=\pm \frac{\pi^2}{4}$. By (\ref{p-eq})
and (\ref{q-eq}),
\[ (1/2)p(t)r''\cos t +(p'(t)\cos t-2p(t)\sin t)r'-r=-4\sin t.
\]
Differentiating the last equation, we get
\begin{eqnarray}
 &[\frac12p(t)\cos t] (r')''+[\frac32p'(t)\cos t-\frac52 p(t)\sin t]
(r')'\nonumber\\
 &+[p''(t)\cos t -3p'(t)\sin t - 2p(t)\cos t -1](r')=-4\cos t.\nonumber
\end{eqnarray}
Using (\ref{p-eq}), the above equation becomes
\begin{eqnarray}                                       \label{r'-eq}
 &[\frac12p(t)\cos t] (r')''+[\frac32p'(t)\cos t-\frac52 p(t)\sin t]
(r')'\nonumber\\
 &+[p'(t)\sin t + 2p(t)\cos t -3](r')=-4\cos t.
\end{eqnarray}
The coefficient of $(r')$ in (\ref{r'-eq}) is negative, for
$p'(t)\sin t +2p\cos t -3< \frac12+ \frac{16}{3\pi} -3 < 0$. This
fact and the negativity of the righthand side of (\ref{r'-eq}) in
$(-\frac{\pi}{2},\frac{\pi}{2})$ imply that $r'$ cannot achieve
its non-positive minimum on $[-\frac{\pi}{2},\frac{\pi}{2}]$ at a
point in $(-\frac{\pi}{2}, \frac{\pi}{2})$.  Now
\begin{eqnarray}\nonumber
  & &\lim_{t\rightarrow \frac{\pi}{2}^-} r'(t) \nonumber\\
  & =&\lim_{t\rightarrow \frac{\pi}{2}^-}s(t)\cos^2 t
  /(\frac{4}{\pi}\cos t +\frac{4}{\pi}t\sin t- \sin^2 t
-1)^2\nonumber\\
 &=&\lim_{t\rightarrow \frac{\pi}{2}^-}[s(t)/ \cos^4]
  /[(\frac{4}{\pi}\cos t +\frac{4}{\pi}t\sin t- \sin^2 t
-1)/ \cos^3 t]^2\nonumber\\
 &=&\lim_{t\rightarrow \frac{\pi}{2}^-}[ s(t)/\cos^4 t]
  /[\frac12\eta'(t)]^2\nonumber\\
  &=&(\frac{4}{3\pi}-\frac{\pi}{12})/(\frac{4}{3\pi})^2\nonumber\\
  &>&0,\nonumber
\end{eqnarray}
where
\begin{eqnarray}\nonumber
s(t)&=&-\frac{4}{\pi}t^2-t^2 \cos t +\frac{12}{\pi}\cos^2 t
  +\frac{8}{\pi}t\sin t\cos t \nonumber\\
  & &-\cos t\sin^2t
 +(\frac{\pi^2}{4}-3)\cos t - \pi + 4t\sin t.\nonumber
\end{eqnarray}
Therefore $r'(t)> 0$ and $r$ is an increasing function on
$[-\frac{\pi}{2},\frac{\pi}{2}]$.
\end{proof}

%
%

Department of mathematics, Utah Valley State College, Orem, Utah
84058

\textit {E-mail address}: \texttt{lingju@uvsc.edu}
\end{document}